\documentclass[10pt,twoside]{article}
\usepackage{graphicx}
\usepackage{amsmath}
\usepackage{amsfonts}
\usepackage{amssymb}
\usepackage{Latex-document}
\newcommand{\ov}[1]{\overline{\vphantom{T}#1}}
\newcommand{\Q}{{\mathbb Q}}
\newcommand{\Z}{{\mathbb Z}}

\newcommand{\eps}{\varepsilon}
\DeclareMathOperator{\NO}{{\cal N}}
\DeclareMathOperator{\N}{{\NO\!}}
\newcommand{\Np}{{\N\p}}
\DeclareMathOperator{\ind}{ind}
\DeclareMathOperator{\SL}{SL}
\DeclareMathOperator{\Sym}{Sym}
\DeclareMathOperator{\rk}{rk}
\newcommand{\gd}{{\mathfrak d}}
\newcommand{\p}{{\mathfrak p}}
\newcommand{\q}{{\mathfrak q}}
\newcommand{\al}{\alpha}
\renewcommand{\a}{{\mathfrak a}}

\renewcommand{\c}{{\mathfrak c}}
\newcommand{\leg}[2]{\mbox{$\left(\frac{#1}{#2}\right)$}}

\title{\bf Constructing and Counting \vskip -2mm  Number Fields\vskip 6mm}
\author{H. Cohen\vspace*{-0.5cm}\thanks{Laboratoire A2X,
Institut de Math\'ematiques, Universit\'e Bordeaux I, 351 Cours de
la Lib\'eration, 33405 TALENCE Cedex, France. E-mail:
cohen@math.u-bordeaux.fr}}
\date{\vspace{-8mm}}

\begin{document}

\maketitle

\thispagestyle{first} \setcounter{page}{129}

\begin{abstract}

\vskip 3mm

In this paper we give a survey of recent methods for the
asymptotic and exact enumeration of number fields with given
Galois group of the Galois closure. In particular, the case of
fields of degree up to 4 is now almost completely solved, both in
theory and in practice. The same methods also allow construction
of the corresponding complete tables of number fields with
discriminant up to a given bound.

\vskip 4.5mm


\end{abstract}

\vskip 12mm

\section{Introduction} \label{section 1}\setzero

\vskip-5mm \hspace{5mm}

Let $K$ be a number field considered as a fixed base field, $\ov{K}$ an
algebraic closure of $K$, and $G$ a transitive permutation group on $n$
letters. We consider the set ${\cal F}_{K,n}(G)$ of all extensions $L/K$ of
degree $n$ with $L\subset\ov{K}$ such that the Galois group of the Galois
closure $\tilde L$ of $L/K$ viewed as a permutation group on the set of
embeddings of $L$ into $\tilde L$ is permutation isomorphic to $G$
(i.e, $n/m(G)$ times the number of extensions up to $K$-isomorphism, where
$m(G)$ is the number of $K$-automorphisms of $L$). We write
$$N_{K,n}(G,X)=|\{L\in{\cal F}_n(G),\ |\N(\gd(L/K))|\le X\}|\;,$$
where $\gd(L/K)$ denotes the relative ideal discriminant and $\N$ the absolute
norm. The aim of this paper is to give a survey of new methods, results, and
conjectures on asymptotic and exact values of this quantity.
It is usually easy to generalize the results to the case where the behavior
of a \emph{finite} number of places of $K$ in the extension $L/K$ is
specified, for example if $K=\Q$ when the signature $(R_1,R_2)$ of $L$ is
specified, with $R_1+2R_2=n$.

{\bf Remarks.}\begin{enumerate}\item It is often possible to give additional main terms and rather good error
terms instead of asymptotic formulas. However, even in very simple cases such as $G=S_3$, this is not at all easy.
\item The methods which lead to exact values of $N_{K,n}(G,X)$ always lead
to algorithms for computing the corresponding \emph{tables}, evidently only when \linebreak $N_{K,n}(G,X)$ is not
too large in comparison to computer memory, see for example \cite{CoDiOlTQA4} and \cite{CoDiOlTQMC}.
\end{enumerate}

General conjectures on the subject have been made by several authors, for
example in \cite{CohBOOK2}. The most precise are due to G.~Malle (see
\cite{Mal1}, \cite{Mal2}). We need the following definition.

{\bf Definition 1.1. }\it For any element $g\in S_n$ different from the
identity, define the index $\ind(g)$ of $g$ by the formula
$\ind(g)=n-|\text{\rm orbits of }g|$.
We define the index $i(G)$ of a transitive subgroup $G$ of $S_n$ by the
formula $$i(G)=\min_{g\in G,\ g\neq1}\ind(g)\;.$$
\rm

{\bf Examples.}\begin{enumerate}\item The index of a transposition
is equal to $1$, and this is the lowest possible index for a
nonidentity element. Thus $i(S_n)=1$.
\item If $G$ is an Abelian group, and if $\ell$ is the smallest prime divisor
of $|G|$, then $i(G)=|G|(1-1/\ell)$.\end{enumerate}

{\bf Conjecture 1.2. }\it For each number field $K$ and transitive
group $G$ on $n$ letters as above, there exist a strictly positive
integer $b_K(G)$ and a strictly positive constant $c_K(G)$ such
that
$$N_{K,n}(G,X)\sim c_K(G)\,X^{1/i(G)}(\log X)^{b_K(G)-1}\;.$$
\rm

In \cite{Mal2}, Malle gives a precise conjectural value for the constant
$b_K(G)$ which is too complicated to be given here.

{\bf Remarks.}\begin{enumerate}\item This conjecture is completely out of reach since it implies the truth of the
inverse Galois problem for number fields.
\item If true, this conjecture implies that for
any \emph{composite} $n$, the proportion of $S_n$-extensions of $K$ of degree
$n$ among all degree $n$ extensions is strictly less than $1$ (but strictly
positive), contrary to the case of \emph{polynomials}.
\end{enumerate}

The following results give support to the conjecture (see \cite{Bha2},
\cite{CoDiOlD4}, \cite{Dat-Wri}, \cite{Dav-Hei1}, \cite{Dav-Hei2},
\cite{Kab-Yuk}, \cite{Klu-Mal}, \cite{Mak1}, \cite{Wri}, \cite{Yuk}).

{\bf Theorem 1.3. }\it We will say that the above conjecture is true in the
weak sense if there exists $c_K(G)>0$ such that for all $\eps>0$ we have
$$c_K(G)\cdot X^{1/i(G)}<N_{K,n}(G,X)<X^{1/i(G)+\eps}\;.$$
\begin{enumerate}\item (M\"aki, Wright). The conjecture is true for all Abelian
groups $G$.
\item (Davenport-Heilbronn, Datskovsky-Wright). The conjecture is true for
$n=3$ and $G=S_3$.
\item (Cohen-Diaz-Olivier). The conjecture is true for $n=4$ and $G=D_4$.
\item (Bhargava, Yukie). The conjecture is true for $n=4$ and $G=S_4$, in the
weak sense if $K\neq\Q$.
\item (Kl\"uners-Malle). The conjecture is true in the weak sense for all
nilpotent groups.
\item (Kable-Yukie). The conjecture is true in the weak sense for $n=5$ and
$G=S_5$.
\end{enumerate}\rm

The methods used to prove these results are quite diverse. In the case of
Abelian groups $G$, one could think that class field theory gives all the
answers so nothing much would need to be done. This is not at all the case, and
in fact Kummer theory is usually more useful. In addition, Kummer theory
allows us more generally to study solvable groups. We will look at this
method in detail.

Apart from Kummer theory and class field theory, the other methods
have a different origin and come from the classification of
\emph{orders} of degree $n$, interpreted through suitable classes
of \emph{forms}. This can be done at a very clever but still
elementary level when the base field is $\Q$, and includes the
remarkable achievement of M.~Bhargava in 2001 for quartic orders.
Over arbitrary $K$, one needs to use and develop the theory of
prehomogeneous vector spaces, initiated at the end of the 1960's
by Sato and Shintani (see for example \cite{Shi1} and
\cite{Shi2}), and used since with great success by
Datskovsky-Wright, and more recently by Wright-Yukie (see
\cite{Wri-Yuk}), Yukie and Kable-Yukie.

\section{Kummer theory} \label{section 2} \setzero

\vskip-5mm \hspace{5mm}

This method applies only to Abelian, or more generally solvable extensions.

\subsection{Why not class field theory?}\vskip-5mm \hspace{5mm}

It is first important to explain why class field theory, which is supposed
to be a complete theory of Abelian extensions, does not give an answer to
counting questions. Let us take the very simplest example of quadratic
extensions, thus with $G=C_2$. A trivial class-field theoretic argument gives
the exact formula
$$N_{K,2}(C_2,X)=-1+\sum_{\N(\a)\le X}2^{\rk(Cl^+_{\a}(K))}M_K\left(\dfrac{X}{\N(\a)}\right)\;,$$
where $\a$ runs over all integral ideals of $K$ of norm less than or equal to
$X$, $Cl^+_{\a}(K)$ denotes the narrow ray class group modulo $\a$,
$\rk(G)$ denotes the $2$-rank of an Abelian group $G$, and
$M_K(n)$ is the generalization to number fields of the summatory function
$M(n)$ of the M\"obius function.

This formula is completely explicit, the quantities $Cl^+_{\a}(K)$ and the
function $M_K(n)$ are algorithmically computable with reasonable efficiency,
so we can compute $N_{K,2}(C_2,X)$ for reasonably small values of $X$ in
this way. Unfortunately, this formula has two important drawbacks.

The first and essential one is that, if we want to deduce from it asymptotic
information on $N_{K,2}(C_2,X)$, we need to control $\rk(Cl^+_{\a}(K)$, which
can be done, although rather painfully, but we also need to control $M_K(n)$,
which \emph{cannot} be done (recall for instance that the Riemann Hypothesis
can be formulated in terms of this function).

The second drawback is that, even for exact computation it is rather
inefficient, compared to the formula that we obtain from Kummer theory.
Thus, even though Kummer theory is used in a crucial way for the constructions
needed in the proofs of class field theory, it must not be discarded once
this is done since the formula that it gives are much more useful, at least
in our context.

\subsection{Quadratic extensions}\vskip-5mm \hspace{5mm}

As an example, let us see how to treat quadratic extensions using
Kummer theory instead of class field theory. Of course in this
case Kummer theory is trivial since it tells us that quadratic
extensions of $K$ are parameterized by $K^*/{K^*}^2$ minus the
unit class. This is not explicit enough. By writing for any
$\al\in K^*$, $\al\Z_K=\a\q^2$ with $\a$ an integral squarefree
ideal, it is clear that $K^*/{K^*}^2$ is in one-to-one
correspondence with pairs $(\a,\ov{u})$, where $\a$ are integral
squarefree ideals whose ideal class is a square, and $\ov{u}$ is
an element of the so-called Selmer group of $K$, i.e., the group
of elements $u\in K^*$ such that $u\Z_K=\q^2$ for some ideal $\q$,
divided by ${K^*}^2$. We can then introduce the Dirichlet series
$\Phi_{K,2}(C_2,s)=\sum_{L}\N(\gd(L/K))^{-s}$, where the sum is
over quadratic extensions $L/K$ in $\ov{K}$. A number of not
completely trivial combinatorial and number-theoretic computations
(see \cite{CoDiOlD4}) lead to the explicit formula
$$\Phi_{K,2}(C_2,s)=-1+\dfrac{2^{-r_2}}{\zeta_K(2s)}\sum_{\c\mid2}\dfrac{\N(2/\c)}{\N(2/\c)^s}\sum_{\chi}L_K(\chi,s)\;,$$
where $\chi$ runs over all quadratic characters of the ray class
group $Cl_{\c^2}(K)$ and $L_K(\chi,s)$ is the ordinary
Dirichlet-Hecke $L$-function attached to $\chi$.

There are two crucial things to note in this formula.
First of all, the sum on $\c$ is only on integral ideals dividing $2$, so
is finite and very small. Thus, $\Phi_{K,2}(C_2,s)$ is a finite linear
combination of Euler products, and can directly be used much more efficiently
than the formula coming from class field theory to compute $N_{K,2}(C_2,X)$
exactly. For example (but this of course does not need the above machinery) we
obtain $N_{\Q,2}(C_2,10^{25})=6079271018540266286517795$.

Second, since $L_K(\chi,s)$ extends to a meromorphic function in the whole
complex plane with no pole if $\chi$ is not a trivial character, the polar
part of $\Phi_{K,2}(C_2,s)$, which is the only thing that we need for an
asymptotic formula, comes only from the contributions of the trivial
characters, in which case $L_K(\chi,s)$ is equal to $\zeta_K(s)$ times
a finite number of Euler factors. We thus obtain
$$N_{K,2}(C_2,X)\sim\dfrac{1}{2^{r_2}}\dfrac{\zeta_K(1)}{\zeta_K(2)}\,X\;,$$
where $\zeta_K(1)$ is a convenient abuse of notation for the
residue of $\zeta_K(s)$ at $s=1$. Apparently this simple result
was first obtained by Datskovsky-Wright in \cite{Dat-Wri},
although their proof is different.

\subsection{General finite Abelian extensions}\vskip-5mm \hspace{5mm}

The same method can in principle be applied to any finite Abelian group $G$.
I say ``in principle'', because in practice several problems arise. For
the base field $K=\Q$, a complete and explicit solution was given by
M\"aki in \cite{Mak1}. For a general base field, a solution has been given
by Wright in \cite{Wri}, but the problem with his solution is that the
constant $c_K(G)$, although given as a product of local
contributions, cannot be computed explicitly without a considerable
amount of additional work. It is always a finite linear combination
of Euler products.

In joint work with F.~Diaz y Diaz and M.~Olivier, using Kummer theory
in a manner analogous but much more sophisticated than the case of
quadratic extensions, we have computed completely explicitly the constants
$c_K(G)$ for $G=C_\ell$ the cyclic group of prime order $\ell$, for
$G=C_4$ and for $G=V_4=C_2\times C_2$. Although our papers are perhaps
slightly too discursive, to give an idea the total number of pages for
these three results exceeds 100. We refer to \cite{CoDiOlCLCRAS},
\cite{CoDiOlCY}, \cite{CoDiOlCL}, \cite{CoDiOlC4}, \cite{CoDiOlV4} for the
detailed proofs, and to \cite{CoDiOlSURVEY} and \cite{CoDiOlCOUNTEXP} for
surveys and tables of results. We mention here the simplest one, for $G=V_4$.
We have
$$N_{K,4}(V_4,X)\sim c_K(V_4)\,X^{1/2}\log^2X\text{\quad with}$$
\begin{align*}c_K(V_4)&=\dfrac{1}{48\cdot 4^{r_2}}\zeta_K(1)^3\prod_{\p}\left(1+\dfrac{3}{\Np}\right)\left(1-\dfrac{1}{\Np}\right)^3\\
&\phantom{=}\prod_{\p\mid2\Z_K}\dfrac{1+\dfrac{4}{\Np}+\dfrac{2}{\Np^2}+\dfrac{1}{\Np^3}-\dfrac{(1-1/\Np^2)e(\p)+(1+1/\Np)^2}{\Np^{e(\p)+1}}}{1+\dfrac{3}{\Np}}\;.\end{align*}
Of course, the main difficulty is to compute correctly the local factor at
$2$.

As usual, we can use our methods to compute very efficiently the $N$ function.
For example, we obtain (see \cite{CohABEL}):
\begin{align*}N_{\Q,3}(C_3,10^{37})&=501310370031289126,\\
N_{\Q,4}(C_4,10^{32})&=1220521363354404,\\
N_{\Q,4}(V_4,10^{36})&=22956815681347605884.
\end{align*}

\subsection{Dihedral {\boldmath $D_4$}-extensions}\vskip-5mm \hspace{5mm}

We can also apply our method to solvable extensions. The case of quartic
$D_4$-extensions, where $D_4$ is the dihedral group of order $8$, is
especially simple and pretty. Such an extension is imprimitive, i.e.,
is a quadratic extension of a quadratic extension. Conversely, imprimitive
quartic extensions are either $D_4$-extensions, or Abelian with Galois
group $C_4$ or $V_4$. These can easily be counted as explained above,
and in any case will not contribute to the main term of the asymptotic
formula, so they can be neglected (or subtracted for exact computations).
Since we have treated completely the case of quadratic extensions, it is
just a matter of showing that we are allowed to sum over quadratic
extensions of the base field to obtain the desired asymptotic formula
(for the exact formula nothing needs to be proved), and this is not
difficult. In this way, we obtain that $N_{K,4}(D_4,X)\sim c_K(D_4)\,X$
for an explicit constant $c_K(D_4)$ (in fact we obtain an error term
$O(X^{3/4}+\eps)$). This result is new even for $K=\Q$, although its
proof not very difficult. In the case $K=\Q$, we have for instance
$$c_{\Q}(D_4)=\dfrac{6}{\pi^2}\sum_{D}\dfrac{2^{-r_2(D)}}{D^2}\dfrac{L(\leg{D}{.},1)}{L(\leg{D}{.},2)}=0.1046520224\dots\;,$$
where the sum is over fundamental discriminants $D$,
$r_2(D)=r_2(\Q(\sqrt{D}))$, and $L(\leg{D}{.},s)$ is the usual Dirichlet
series for the character $\leg{D}{.}$.

{\bf Remark.} In the Abelian case, it is possible to compute the Euler
products which occur to hundreds of decimal places if desired using
almost standard zeta-product expansions, see for example \cite{CohLW}.
Unfortunately, we do not know if it is possible to express $c_{\Q}(D_4)$
as a finite linear combination of Euler products (or at least as a
rapidly convergent infinite series of such), hence we have only been able
to compute 9 or 10 decimal places of this constant. We do not see any
practical way of computing 20 decimals, say.

Our method also allows us to compute $N_{\Q,4}(D_4,X)$ exactly. However,
here a miracle occurs: when $k$ is a quadratic field, in the formula that
we have given above for $\Phi_{k,2}(C_2,s)$ all the quadratic characters
$\chi$ which we need are \emph{genus characters} in the sense of Gauss,
in other words there is a decomposition
$$L_k(\chi,s)=L(\leg{d_1}{.},s)L(\leg{d_2}{.},s)$$
into a product of two suitable ordinary Dirichlet $L$-series.
This gives a very fast method for computing $N_{\Q,4}(D_4,X)$, and
in particular we have been able to compute
$N_{\Q,4}(D_4,10^{17})=10465196820067560$.

We can also count the number of extensions with a given signature. The
method is completely similar, but here not all characters are genus
characters. In fact, it is only necessary to add a single nongenus
character to obtain all the necessary ones, but everything is completely
explicit, and closely related to the \emph{rational quartic reciprocity
law}. I refer to \cite{CohD4REAL} for details.

\subsection{Other solvable extensions}\vskip-5mm \hspace{5mm}

We can also prove some partial results in the case where $G=A_4$ or
$G=S_4$ (of course the results for $S_4$ are superseded by Bhargava's
for $K=\Q$, and by Yukie's for general $K$; still, the method is also
useful for exact computations), see \cite{CoDiOlA4S4}.

In the case of quartic $A_4$ and $S_4$-extensions (or, for that matter, of
cubic $S_3$-extensions), we use the diagram involving the cubic resolvent
(the quadratic one for $S_3$-extensions), also called the Hasse diagram. We
then have a situation which bears some analogies with
the $D_4$ case. The differences are as follows. Instead of having to sum over
quadratic extensions of the base field $K$, we must sum over cubic extensions,
cyclic for $A_4$ and noncyclic for $S_4$. As in the $D_4$-case, we then have to
consider quadratic extensions of these cubic fields, but generated by an
element of square norm. It is possible to go through the exact combinatorial
and arithmetic computation of the corresponding Dirichlet series, the cubic
field being fixed. This in particular uses some amusing local class field
theory. As in the $D_4$ case, we then obtain the Dirichlet generating series
for discriminants of $A_4$ (resp., $S_4$) extensions by summing the series
over the corresponding cubic fields.

Unfortunately, we cannot obtain from this any asymptotic formula. The reason
is different in the $A_4$ and the $S_4$ case. In the $A_4$ case, the
rightmost singularity of the Dirichlet series is at $s=1/2$. Unfortunately,
this is simultaneously the main singularity of each individual Dirichlet
series, and also that of the generating series for cyclic cubic fields.
Thus, although the latter is well understood, it seems difficult (but not
totally out of reach) to paste things together. On the other hand, we can
do two things rigorously in this case. First, we can prove an asymptotic
formula for $A_4$-extensions having a \emph{fixed} cubic resolvent. Tables
show that the formula is very accurate. Second, we can use our formula to
compute $N_{K,4}(S_4,X)$ exactly. For instance, we have computed
$N_{\Q,4}(A_4,10^{16})=218369252$. This computation is much slower than in
the $D_4$-case, because we do not have the miracle of genus characters, and
we must compute the class and unit group of all the cyclic cubic fields.

In the $S_4$ case, the situation is different. The main
singularity of each individual Dirichlet series is still at
$s=1/2$ (because of the square norm condition), and the rightmost
singularity of the generating series for noncyclic cubic fields is
at $s=1$, so the situation looks better (and analogous to the
$D_4$ situation with $s$ replaced by $s/2$). Unfortunately, as
already mentioned we know almost nothing about the generating
series for noncyclic cubic fields, a fortiori with coefficients.
So we cannot go further in the asymptotic analysis. As in the
$A_4$ case, however, we can compute exactly either the number of
$S_4$-extensions corresponding to a fixed cubic resolvent, or even
$N_{K,4}(S_4,X)$ itself. The problem is that here we must compute
class and unit groups of all noncyclic cubic fields of
discriminant up to $X$, while cyclic cubic fields of discriminant
up to $X$ are much rarer, of the order of $X^{1/2}$ instead. We
have thus not been able to go very far and obtained for example
$N_{\Q,4}(S_4,10^7)=6541232$.

\section{Prehomogeneous vector spaces} \label{section 3}
\setzero\vskip-5mm \hspace{5mm }

The other methods for studying $N_{K,n}(G,X)$ are two closely
related methods: one is the use of generalizations of the
Delone-Fadeev map, which applies when $K=\Q$. The other, which can
be considered as a generalization of the first, is the use of the
theory of prehomogeneous vector spaces, initiated by Sato and
Shintani in the 1960's.

\subsection{Orders of small degree}\vskip-5mm \hspace{5mm}

We briefly give a sketch of the first method. We would first like to
classify \emph{quadratic} orders. It is well known that, through their
discriminant, such orders are in one-to-one correspondence with the subset
of nonsquare elements of $\Z$ congruent to $0$ or $1$ modulo $4$, on which
$\SL_1(\Z)$ (the trivial group) acts. Thus, for fixed discriminant, the
orbits are finite (in fact of cardinality $0$ or $1$).
For \emph{maximal} orders, we need to add local arithmetic conditions at each
prime $p$, which are easy for $p>2$, and slightly more complicated for $p=2$.

We do the same for small higher degrees. For \emph{cubic} orders,
the classification is due to Davenport-Heilbronn (see
\cite{Dav-Hei1}, \cite{Dav-Hei2}). These orders are in one-to-one
correspondence with a certain subset of $\Sym^3(\Z^2)$, i.e.,
binary cubic forms, on which $\SL_2(\Z)$ acts. Since once again
the difference in ``dimensions'' is $4-3=1$, for fixed
discriminant the orbits are finite, at least generically. For
maximal orders, we again need to add local arithmetic conditions
at each prime $p$. These are easy to obtain for $p>3$, but are a
little more complicated for $p=2$ and $p=3$. An alternate way of
explaining this is to say that a cubic order can be given by a
\emph{nonmonic} cubic equation, which is almost canonical if
representatives are suitably chosen.

For \emph{quartic} orders, the classification is due to M.~Bhargava in 2001,
who showed in complete detail how to generalize the above. These orders are
now in one-to-one correspondence with a certain subset of
$\Z^2\otimes\Sym^2(\Z^3)$, i.e., pairs of ternary quadratic forms, on
which $\SL_2(\Z)\times\SL_3(\Z)$ acts. Once again the difference in
``dimensions'' is $2\times6-(3+8)=1$, so for fixed discriminant the orbits
are finite, at least generically. For maximal orders, we again need to add
local arithmetic conditions at each prime $p$, which Bhargava finds after
some computation. An alternate way of explaining this is to say that a
quartic order can be given by the intersection of two conics in the projective
plane, the pencil of conics being almost canonical if representatives are
suitably chosen.

For \emph{quintic} orders, only part of the work has been done, by
Bhargava and Kable-Yukie in 2002. These are in one-to-one
correspondence with a certain subset of
$\Z^4\otimes\Lambda^2(\Z^5)$, i.e., quadruples of alternating
forms in $5$ variables, on which $\SL_4(\Z)\times\SL_5(\Z)$ acts.
Once again the difference in ``dimensions'' is
$4\times10-(15+24)=1$, so for fixed discriminant the orbits are
finite, at least generically. The computation of the local
arithmetic conditions, as well as the justification for the
process of point counting near the cusps of the fundamental domain
has however not yet been completed.

Since prehomogeneous vector spaces have been completely classified, this
theory does not seem to be able to apply to higher degree orders, at least
directly.

\subsection{Results}\vskip-5mm \hspace{5mm}

Using the above methods, and generalizations to arbitrary base fields, the
following results have been obtained on the function $N_{K,n}(G,X)$ (many
other deep and important results have also been obtained, but we fix our
attention to this function). It is important to note that they
seem out of reach using more classical methods such as Kummer theory or
class field theory mentioned earlier.

{\bf Theorem 3.1. }\it Let $K$ be a number field of signature $(r_1,r_2)$, and
as above write $\zeta_K(1)$ for the residue of the Dedekind
zeta function of $K$ at $s=1$.
\begin{enumerate}\item (Davenport-Heilbronn \cite{Dav-Hei1}, \cite{Dav-Hei2}).
We have $N_{\Q,3}(S_3,X)\sim c_{\Q}(S_3)\,X$ with
$$c_{\Q}(S_3)=\dfrac{1}{\zeta(3)}\;.$$
\item (Datskovsky-Wright \cite{Dat-Wri}). We have
$N_{K,3}(S_3,X)\sim c_K(S_3)\,X$ with
$$c_K(S_3)= \left(\dfrac23\right)^{r_1-1}\left(\dfrac16\right)^{r_2}\dfrac{\zeta_K(1)}{\zeta_K(3)}\;.$$
\item (Bhargava \cite{Bha1}, \cite{Bha2}). We have
$N_{\Q,4}(S_4,X)\sim c_{\Q}(S_4)\,X$ with
$$c_{\Q}(S_4)=\dfrac{5}{6}\prod_p\left(1+\dfrac1{p^2}-\dfrac1{p^3}-\dfrac1{p^4}\right)\;.$$
\item (Yukie \cite{Yuk}). There exist two strictly positive constants $c_1(K)$
and $c_2(K)$ such that
$$c_1\,X<N_{K,4}(S_4,X)<c_2\,X\log^2(X)\;.$$
Under some very plausible convergence assumptions we should have in fact
$N_{K,4}(S_4,X)\sim c_K(S_4)\,X$ with
$$c_K(S_4)=2\left(\dfrac{5}{12}\right)^{r_1}\left(\dfrac1{24}\right)^{r_2}\prod_{\p}\left(1+\dfrac1{\N\p^2}-\dfrac1{\N\p^3}-\dfrac1{\N\p^4}\right)\;.$$
\item (Kable-Yukie \cite{Kab-Yuk}). There exists a strictly positive constant
$c_1$ such that for all $\eps>0$ we have
$$c_1\,X<N_{\Q,5}(S_5,X)<X^{1+\eps}\;.$$
\end{enumerate}\rm

{\bf Remark.} It should also be emphasized that, although the above methods give important and deep results on
$N_{K,n}(G,X)$ for certain groups $G$, they shed almost no light on the possible analytic continuation of the
corresponding Dirichlet series of which $N_{K,n}(G,X)$ is the counting function. For example, in the simplest case
where $K=\Q$, $n=3$, and $G=S_3$, for which the result dates back to Davenport-Heilbronn, no one knows how to give
an analytic continuation of the Dirichlet series $\sum_{L}|d(L)|^{-s}$ even to $\Re(s)=1$ (the sum being over
cubic fields in $\ov{\Q}$ and $d(L)$ being the absolute discriminant of $L$).

\label{lastpage}

\end{document}